\documentclass{article}

\usepackage{amsfonts}
\usepackage{amssymb}
\usepackage{amsmath}
\usepackage{amsthm}
\usepackage{latexsym}
\usepackage{graphicx}

\newtheorem{thm}{Theorem}

\newtheorem{claim}[thm]{Claim}
\numberwithin{equation}{section}
\theoremstyle{definition}

\usepackage{color}

\newcommand{\lp}{\left(}
\newcommand{\rp}{\right)}

\newcommand{\pd}[2]{\frac{\partial #1}{\partial #2}}
\newcommand{\spd}[2]{\frac{\partial^2 #1}{\partial {#2}^2}}

\definecolor{dred}{rgb}{.9,0,.1}
\newcommand{\mR}{\mathbb{R}}

\newcommand{\R}{\mathbb{R}}

\newcommand{\C}{\mathbb{C}}

\newcommand{\aA}{\mathcal A}
\newcommand{\bB}{\mathcal B}
\newcommand{\cC}{\mathcal C}
\newcommand{\tT}{\mathcal T}

\newcommand{\be}{\begin{equation}}
\newcommand{\ee}{\end{equation}}
\newcommand{\bee}{\begin{equation*}}
\newcommand{\eee}{\end{equation*}}
\newcommand{\bea}{\begin{eqnarray}}
\newcommand{\eea}{\end{eqnarray}}
\newcommand{\bs}{\begin{split}}
\newcommand{\es}{\end{split}}

\newcommand{\imagpart}{{\mathcal Im}}

\begin{document}

\title{{\bf Instability and bifurcation in a trend depending price formation model}}

 \author{Mar\'ia del Mar Gonz\'alez \thanks{Partially supported by grants
MINECO MTM2011-27739-C04-01 and GENCAT 2009SGR-345.}\\Univ. Polit\`ecnica de Catalunya  \and Maria Pia Gualdani \thanks{Supported by the NSF Grant DMS-1109682.}\\George Washington University  \and Joan Sol\`{a}-Morales  \thanks{Partially supported by grants MINECO MTM2011-27739-C04-01 and GENCAT 2009SGR-345.}\\Univ. Polit\`ecnica de Catalunya}

\date{}

\maketitle


\begin{abstract}
A well-known model due to J.-M. Lasry and P.L. Lions that presents the evolution of prices in a market as the evolution of a free boundary in a diffusion equation is suggested to be modified in order to show instabilities for some values of the parameters. This loss of stability is associated to the appearance of new types of solutions, namely periodic solutions, due to a Hopf bifurcation and representing price oscillations, and traveling waves, that represent either inflationary or deflationary behavior.
\end{abstract}


2000 {\bf MSC:} 35R35, 35K15, 91B42, 91B26

%

\section{Introduction and summary of results}

A mathematical model for the time evolution of a price in some trading markets, where the actual price is the location of a free boundary of a nonlinear diffusion problem, was proposed by J.-M. Lasry and P.-L. Lions \cite{Lasry-Lions} in 2007. In that model, the population in the market is divided into two groups, namely the buyers ($B$) and vendors ($V$), and described by two time-dependent densities $f_B(x,t)$ and $f_V(x,t)$ respectively, where $x$ has the dimensions of a price and $f_B(x,y)$ and $f_V(x,t)$ respectively mean the amount of buyers or sellers that would accept a transaction at price $x$. The value of $x$ at which the actual transaction takes place is denoted by $p(t)$, and, by definition, $f_B(x,t)=0$ for $x>p(t)$ and $f_V(x,t)=0$ for $x<p(t)$.

The time evolution of these densities obeys first to a diffusion law, coming from the idea that both buyers and vendors change their minds on the desired prices following a Gaussian random process with variance $\sigma^2$. And second, when some buyers agree on prices higher than $p(t)$ and/or some vendors accept prices lower that $p(t)$, then new sales happen, which produces a change in the location of the free boundary $x=p(t)$. Then, the former buyers and sellers leave the market with a flux denoted by $\Lambda(t)$ and immediately re-enter as new sellers and buyers, respectively, at values of  $x=p(t)+a$ and $x=p(t)-a$, where $a>0$ is the so-called transaction cost.

These ideas are reflected in the following equations:

\begin{equation}
\label{equations}\left\{ \begin{array}{ll}
&\dfrac{\partial f_B}{\partial t} -\dfrac{\sigma^2}{2}\dfrac{\partial^2f_B}{\partial x^2}=\Lambda(t)\delta_{p(t)-a}
\quad \hbox{if}\: x< p(t), \; t>0,\\\\
&\dfrac{\partial f_V}{\partial t} -\dfrac{\sigma^2}{2}\dfrac{\partial^2f_V}{\partial x^2}=\Lambda(t)\delta_{p(t)+a}
\quad \hbox{if } x>p(t), \; t>0,\\\\
&\Lambda(t)=-\dfrac{\sigma^2}{2}\dfrac{\partial f_B}{\partial x}(p(t),t)=\dfrac{\sigma^2}{2}\dfrac{\partial f_V}{\partial x}(p(t),t).
\end{array}\right.
\end{equation}
with $f_B(p(t),t)=f_V(p(t),t)=0$, $f_B(x,t)>0$ for $x<p(t)$, $f_V(x,t)>0$ for $x>p(t)$ and one also may suppose that both $f_B(x,t)$ and $f_V(x,t)$ remain bounded as $|x|\to\infty$. We used the symbol $\delta_{x_0}$ to mean a Dirac delta function centered at $x=x_0$. These equations have to be complemented with suitable initial conditions.

The equations (\ref{equations}) can be reduced to a single one with the new unknown $f(x,t)=f_B(x,t)-f_V(x,t)$, namely
\begin{equation}
\label{single}
f_t-\dfrac{\sigma^2}{2}f_{xx}=-\dfrac{\sigma^2}{2}f_x(p(t),t)\left(\delta_{p(t)-a}-\delta_{p(t)+a}\right)
\end{equation}
with the additional condition $f(p(t),t)=0$ and $f(x,t)>0$ for $x<p(t)$ and $f(x,t)<0$ for $x>p(t)$.

Taking derivatives with respect to time in the expression $f(p(t),t)=0$ and using that $f_t(p(t),t)=\frac{\sigma^2}{2}f_{xx}(p(t),t)$ we see that we can use the alternative expression for $p'(t)$
\begin{equation}
\label{pprime}
p'(t)= -\frac{\sigma^2}{2}f_{xx}(p(t),t)/f_x(p(t),t)
\end{equation}
Observe in (\ref{pprime}) that the law that governs the time evolution of the free boundary is nonlinear, and depends on the second derivative of the unknown.

The same equations have also been considered in a bounded domain $a_0<x<a_1$ and then the conditions at infinity have been replaced by homogeneous Neumann boundary conditions at $x=a_0,a_1$. In this Neumann case one easily sees that the two total masses $$\int_{a_0}^{p(t)}f_B(x,t)\ dx\ \ \ \hbox{  and  }\ \ \ \int_{p(t)}^{a_1}f_V(x,t)\ dx$$ remain constant in time.

Both the dynamics of the problem on the whole line and that of the Neumann problem are actually reasonably well understood by the work done in the last few years by several authors, like M.d.M. Gonz\'{a}lez, M.P. Gualdani, L. Chayes, I. Kim, P.A. Markowich, N. Matevoysan, J.-F. Pietchmann, M.-T. Wolfram and L.A. Caffarelli (\cite{Gonzalez-Gualdani:symmetric,Chayes-Gonzalez-Gualdani-Kim,Markowich:price-formation,
Gonzalez-Gualdani:asymptotics,Gonzalez-Gualdani:spaces,Caffarelli-Markowich-Pietschmann,
Caffarelli-Markowich-Wolfram}). Essentially, it has been shown that every solution approaches a single equilibrium when time tends to infinity. This implies, in particular, that there are no solutions of the form of a traveling wave, something that can also be checked directly. A modification of this model that considers possible deals outside the fixed price, with a dynamics of a  Boltzman-type collision, has recently been studied by M. Burger, L.A. Caffarelli P.A. Markowich and M.-T. Wolfram in \cite{Burger-Caffarelli-Markowich-Wolfram}.

It is easy to see that the equilibrium solutions of (\ref{single}) are the piecewise linear functions of the form $f(x)=w^\rho(x-p_0)$, where
\begin{equation}
\label{equilibrium}
w^{\rho}(x)= \begin{cases}
-\rho x/a& \text{ if $|x|\le a$},\\
\rho& \text{if $x<-a$},\\
-\rho& \text{if $x>a$},
\end{cases}
\end{equation}
and $\rho,p_0$ are parameters, $\rho>0$, $p_0\in\R$.

According to this model, markets always stabilize. One could argue that this is not the case in real life. Our aim in this paper is to present modifications of the model (\ref{equations}) that without changing the equilibrium solutions make them to become unstable, at least in some cases.

An inspiring example of destabilization of an equilibrium in a diffusion equation has been that of P. Guidotti and S. Merino \cite{Guidotti-Merino}. In that example a kind of heat regulation mechanism produces oscillations in the temperature when a parameter becomes sufficiently large. These oscillations appear as a consequence of a Hopf bifurcation. Following in part this idea we propose the following new model to replace (\ref{equations}):

\begin{equation}
\label{new}\left\{ \begin{array}{ll}
&\dfrac{\partial f_B}{\partial t} -\dfrac{\sigma^2}{2}\dfrac{\partial^2f_B}{\partial x^2}=\left(\Lambda(t)-Rp'(t)\right)\delta_{p(t)-a}
\quad \hbox{if}\: x< p(t), \; t>0,\\\\
&\dfrac{\partial f_V}{\partial t} -\dfrac{\sigma^2}{2}\dfrac{\partial^2f_V}{\partial x^2}=\left(\Lambda(t)+Rp'(t)\right)\delta_{p(t)+a}
\quad \hbox{if } x>p(t), \; t>0,\\\\
&\Lambda(t)=-\dfrac{\sigma^2}{2}\dfrac{\partial f_B}{\partial x}(p(t),t)=\dfrac{\sigma^2}{2}\dfrac{\partial f_V}{\partial x}(p(t),t).
\end{array}\right.
\end{equation}
with $f_B(p(t),t)=f_V(p(t),t)=0$, $f_B(x,t)>0$ for $x<p(t)$, $f_V(x,t)>0$ for $x>p(t)$ and one also supposes that both $f_B(x,t)$ and $f_V(x,t)$ remain bounded as $|x|\to\infty$.
In this model, the evolution of the prices affects the behavior of the market, and this is why we say this model may be called trend dependent.

A reaction term of the form $Rp'(t)\delta_{p(t)\pm a}$ also recalls the nonlinear part of the equation studied in \cite{CarrilloEtAl} after a suitable transformation.

Here $R$ is a parameter. When $R=0$ one recovers the original model (\ref{equations}). The meaning of $R>0$ is that if the prices grow ($p'(t)>0$) then some buyers leave the market and at the same time some people outside the market (perhaps these previous buyers) enter as vendors. This is somehow the naive idea that when the prices are high it is time to sell, not to buy. And the contrary if $p'(t)<0$. This is a kind of non-local regulatory behavior of the type of that of \cite{Guidotti-Merino}. Of course one could argue that there are other places, and not only $x=p(t)\pm a$, to leave or enter the market, but we think of this as the possibility that makes the simplest model.

But the case $R<0$, being right the contrary, is also meaningful. If $p'(t)>0$ (prices growing) it implies that some people outside the market enter into the game as buyers, perhaps because they feel that the prices may keep growing for some time, so it is a good moment to buy. And in the same situation ($R<0$ and $p'(t)>0$) some vendors leave the market, perhaps also expecting the prices to keep growing and re-enter into the market as vendors when the prices become higher. And the contrary, if $R<0$ and $p'(t)<0$.

In summary, somehow, $R>0$ means conservative market, while a more aggressive investment is represented by $R<0$. Without being precise, one can roughly say that our results with $R>0$ will lead to oscillations and, on the contrary, $R<0$ will lead to traveling waves, both inflationary and deflationary.

The present paper should be merely seen as a mathematical discussion of some simple ways to destabilize the original model (\ref{equations}) and we do not make attempt to have a rigorous financial discussion. A more precise formulation would include a discussion on the possible reasonable boundary conditions, the review of the ideas under the transaction cost $a>0$, and possible alternatives for the nonlinear reaction term.

One can write (\ref{new}) as a single equation for $f=f_B-f_V$ as it has been done in (\ref{single}),
\begin{equation}
\label{newsingle}
f_t-\dfrac{\sigma^2}{2}f_{xx}=-\dfrac{\sigma^2}{2}f_x(p(t),t)\left(\delta_{p(t)-a}-\delta_{p(t)+a}\right)-Rp'(t)\left(\delta_{p(t)-a}+\delta_{p(t)+a}\right).
\end{equation}

Let us discuss briefly possible alternative nonlinearities. Numerical experiments show that when the equilibrium solutions become unstable sometimes happens that the function $f(x,t)$ loses the right sign near $x=p(t)-a$ and $x=p(t)+a$, becoming, respectively, negative and positive and then physically nonsense. To avoid this behavior one can substitute the term $-Rp'(t)\left(\delta_{p(t)-a}+\delta_{p(t)+a}\right)$ in (\ref{newsingle}) by the new similar expression $-Rp'(t)\left(f(p(t)-a,t)\delta_{p(t)-a}-f(p(t)+a,t)\delta_{p(t)+a}\right)$ or, with more generality, by
\begin{equation}
\label{nonlinearity}
-Rp'(t)\left(\phi(f(p(t)-a,t))\delta_{p(t)-a}-\phi(f(p(t)+a,t))\delta_{p(t)+a}\right),
\end{equation}
where $\phi(r)$ satisfies $\phi(-r)=-\phi(r)$ and $\phi(r)>0$ when $r>0$. Natural choices for $\phi$ could be $\phi(r)=\hbox{sign}(r)$, and one recovers (\ref{newsingle}), $\phi(r)=r$ or $\phi(r)=\tanh(r)$. This last choice is inspired by \cite{Guidotti-Merino}, and we can say it is the nonlinearity that gives the clearest numerical results.

Another factor to take into account in the choice of $\phi(r)$ is the possible invariance of the problem under multiplication by a positive constant. More precisely, in the original problem (\ref{single}) if $f(x,t)$ is a solution, then also $\mu f(x,t)$ is a solution, for $\mu>0$. This property is preserved, and with the same value of $R$, for the nonlinearity $\phi(r)=r$. For the other choices of the nonlinearity one should not forget that this property is lost, which introduces a new ingredient into the modeling discussion.

Our numerical experiments with these new models show that the equilibrium solutions become unstable for large values of $R$, both positive and negative. Also, bifurcated solutions appear, in the form of periodic oscillations for $R>0$, consequence of a Hopf bifurcation, and solutions in the form of traveling waves when $R<0$.

In the present paper we provide analytical arguments that support this instability and these bifurcations. But we accept from the beginning that some of these analytical arguments are not complete, and deserve further study, with more modeling discussion and more numerical simulations. Instability and bifurcation in parabolic equations are well-known, but delicate issues (see Henry \cite{Henry:parabolic}, Chow-Hale \cite{Chow-Hale:bifurcation-theory}, Haragus-Iooss \cite{Haragus-Iooss}), whose analysis relies, for example, in appropriate choices of the function spaces, which is particularly delicate in an unbounded domain. The presence of essential spectrum of the second derivative operator up to $\lambda =0$ in the complex plane makes bifurcation analysis more complicated, and in our problem this appears together with multiple eigenvalues also at $\lambda=0$. The reduction to a bounded domain would simplify the first problem, but not the second. And by the way it would make the existence of traveling wave solutions something less clear. Let us present now our results.

Without loss of generality we suppose $\sigma^2/2=1$ and $a=1$. We rewrite equation (\ref{newsingle}) for a new unknown $w(x,t)=f(x+p(t),t)$ that locates the free boundary always in $x=0$:
\begin{equation}\label{eqnonlin1}
\left\{\begin{split}
&w_t =w_{xx} +p'(t)w_x(x,t)-w_x(0,t)\left[\delta_{-1}-\delta_{1}\right]-Rp'(t)\left[\phi(w(-1,t))\delta_{-1}
-\phi(w(1,t))\delta_{1}\right],\\
&w(0,t)=0, \ p'(t)=-w_{xx}(0,t)/w_x(0,t).
\end{split}\right.
\end{equation}

We also concentrate our attention on the stability of and the bifurcation from the equilibrium $w^1(x)$ defined by (\ref{equilibrium}) with $\rho=1$. Other equilibrium solutions can be obtained from $w^1(x)$ by a scale factor, that will enter into the equation either as a change in $R$ or in $\phi(r)$. To simplify again, and without loss of generality, we will consider only nonlinearities of the form (\ref{nonlinearity}) such that $\phi(1)=-\phi(-1)=1$. In fact we have in mind the three principal cases:
\begin{equation}
\label{nonlinearities}\phi(r)=\left\{ \begin{split}
&\phi_1(r)=\hbox{ sign}(r), \ \ \hbox{or}\\
&\phi_2(r)=r, \ \ \hbox{or}\\
&\phi_3(r)=\tanh(r)/\tanh(1).
\end{split}\right.
\end{equation}

The rest of the paper is organized in two more Sections. Section 2 is devoted to an analysis of the stability and instability of the equilibrium $w^1(x)$ defined above by looking at the eigenvalues of the linear part around it. It is clear that, since the functional framework is not completely fixed, one cannot speak about points of the spectrum that are not eigenvalues. Even the linearization itself has to be understood in a very naive way. We will calculate the eigenvalues accepting as a definition that the corresponding eigenfunction is globally bounded. Without trying to be very rigorous, it is known from other cases that when one restricts oneself to eigenfunctions that vanish at infinity these eigenvalues do not longer exist, but at least in some cases they remain in the spectrum as parts of a continuous spectrum.

One should note, by the way, that the linearization around $w^1(x)$ turns out to be the same disregarding the choice of the nonlinearity $\phi(r)$, as long as $\phi(\pm 1)=\pm 1$.

With these considerations in mind, we can summarize the results of Section 2 in the following

\begin{thm}\label{thm1}
The eigenvalues $\lambda$ of the linear part of (\ref{eqnonlin1}) around $w=w^1$, given by
\begin{equation}\label{linearized}
 L_1g:=g_{xx}- g_{xx}(0,t)\mathcal X_{(-1,1)}-g_x(0,t)[\delta_{-1} -\delta_{1}] + Rg_{xx}(0,t)  \left[\delta_{-1}+\delta_{1}\right],
 \end{equation}
 satisfy $$\Re(\lambda)\le 0 \hbox{ for all } -1\le R\le R_0,$$ where $R_0\simeq 9.36\ldots=(1-e^{a_0}\cos(a_0))/a_0$ and $a_0$ is the root of $\cos(a)-\sin(a)=e^{-a}$ near $a\simeq 3.94\dots$.
For $R=-1$ a real eigenvalue becomes positive, and for $R=R_0$ a pair of simple complex conjugated nonzero eigenvalues cross the imaginary axis from left to right.
\end{thm}

In (\ref{linearized}) $\mathcal X_{(-1,1)}(x)$ means the characteristic function of the interval $(-1,1)$.

Section 3 deals with solutions bifurcating from $w^1(x)$ at these values $R=-1$ and $R=R_0$. Roughly speaking, for $R=-1$ the bifurcated solutions are traveling waves, moving both right and left, and for $R=R_0$ an Andronov-Hopf type bifurcation occurs, giving rise to periodic oscillations. This is summarized in the following two results:

\begin{thm}\label{thm2}

For certain ranges of $R<0$ a global two-parameter family of traveling waves exist, that bifurcates from the one-dimensional family of equilibria. This family bifurcates from $w^1(x)$ at $R=-1$. These waves appear in pairs, with velocity $c>0$, meaning an inflationary solution, and with velocity $-c$, that is deflationary. For the three nonlinearities of (\ref{nonlinearities}) they are described more explicitly. For $\phi=\phi_1$ they occupy the whole range $R\in(-\infty, 0)$, for $\phi=\phi_2$ this continuum of solutions lies entirely in $R=-1$, while for $\phi=\phi_3$ it occupies the range $-\infty<R<-\tanh(1)$.
\end{thm}

\begin{claim}\label{claim}
We have numerically observed that at $R=R_0$ a family of periodic solutions does appear near $w^1$. These numerical simulations show at least that for $\phi=\phi_3$ the bifurcation is supercritical, and stable oscillations seem to persist for all $R>R_0$.
\end{claim}

\section{Analysis of the linear problem}

We consider the equilibrium solution $w^1(x)$ defined in (\ref{equilibrium}).
To linearize (\ref{eqnonlin1}) around it we write $w=w^1+\varepsilon g$, substitute in (\ref{eqnonlin1}), differentiate with respect to $\varepsilon$ and set $\varepsilon=0$. That yields to
\begin{equation}\label{eqlin1}
g_t =L_1g,
\end{equation}
for $L_1$ defined in \eqref{linearized}, as long as $\phi(1)=-\phi(-1)=1$, where $\mathcal X_{(-1,1)}(x)=w^0_x(x)/w^0_x(0)$. Since the function $w(x)=w^0(x)+\varepsilon g(x)$ has to satisfy $w(0)=0$ one has also to impose $g(0)=0$ to the perturbation part.

To prove Theorem \ref{thm1} we have to find all the pairs $(g(x),\lambda)$ such that $L_1g=\lambda g$  and $g(x)$ is a nonzero bounded function defined in $-\infty<x<\infty$ with $g(0)=0$. It is easy to see that the equation $L_1g=\lambda g$ decouples if one splits $g(x)$ into its even and odd part. Thus, we write $g = g_1+g_2$ where $g_1$ is even and $g_2$ is odd and we substitute $g_1 + g_2$ into the eigenvalue equation. One gets
\begin{align*}
{g_1}_{xx}  -{g_1}_{xx}(0)\mathcal{X}_{(-1,1)}-R{g_1}_{xx}(0)\left[\delta_{-1}+\delta_{1}\right]=\lambda g_1, \\
{g_2}_{xx}  -{g_2}_{x}(0)\left[\delta_{-1}-\delta_{1}\right]=\lambda g_2,
\end{align*}
equations that need now to be solved for globally bounded solutions in $(-\infty,\infty)$ together with the conditions $g_1(0)=g_{1x}(0)=0$ in the first case and $g_{2}(0)=0$ in the second case. At a first sight we already see that since the equation for the odd part does not depend on $R$ one should only expect unstable modes among even functions of $x$.

Now we proceed to calculate the eigenvalues. We start with the odd eigenfunctions. It is easy to see that for $\lambda=0$ there exists a single bounded
eigenfunction with the same form as the equilibria (\ref{equilibrium}). It is the natural zero eigenvalue whose eigenfunction is tangent to the curve of equilibria.

So we may assume that $\lambda = \alpha^2$ with $\alpha\in \C\backslash\{0\}$ and split $g_2$ as
\bee
g_2(x)=
\left\{
\begin{split}
&h_1(x)=c_1 e^{\alpha x}+c_2 e^{-\alpha x}&\quad\mbox{ in }(-1,1),\\
&h_2(x)=d_1 e^{\alpha x}+d_2 e^{-\alpha x} &\quad\mbox{ in }(1,\infty),\\
&h_3(x)=-d_2 e^{\alpha x}-d_1e^{-\alpha x}&\quad\mbox{ in }(-\infty,-1)\\
\end{split}\right.\eee
for some constants $c_1,c_2,d_1$ and $d_2$. Since we seek odd and bounded eigenfunctions, to solve the above problem in the whole line is equivalent to solve it in the half line with the additional Dirichlet condition $g_2(0)=0$. Hence $h_1$, $h_2$ and $h_3$ have to satisfy the following matching and jump conditions:
\bee
\left\{
\begin{split}
&h_1(0)=0,\\
&{h_2}_x(1^+) - {h_1}_x(1^-) = - {h_1}_x(0), \\
&h_1(1) = h_2(1).\\
\end{split}\right.\eee
These matching conditions imply that the coefficients $d_1$, $d_2$, $c_1$, $c_2$ and $\alpha$ satisfy the following system
\bee
\left\{
\begin{split}
&c_1+c_2 =0,\\
&d_1+d_2e^{-2\alpha} = c_1 + c_2 e^{-2\alpha}, \\
&d_1e^{\alpha} - d_2 e^{-\alpha}  - c_1 e^{\alpha} + c_2 e^{-\alpha} = -c_1+c_2.
\end{split}\right.\eee
We distinguish now two cases: $\Re(\alpha )>0$ (i.e. $d_1=0$) and $\Re (\alpha)=0$. The case $\Re(\alpha )< 0$ needs not to be considered if we have already considered $\Re(\alpha )> 0$ since we are only interested in $\lambda=\alpha^2$.
 A simple calculation shows that the possible nontrivial solutions can only exist with $\Re(\alpha) = 0$. Hence $\alpha = i b$ with $b\in (0,\infty)$, since as before the cases $b<0$ are also considered in $b>0$. Consequently any eigenvalue $\lambda$ is real and negative.
In the particular cases that $\alpha_k=2k\pi i$, $k\in\mathbb Z\backslash\{0\}$, the eigenfunction is
\begin{equation*}
     g_2^k(x)=
     \left\{\begin{split}
     c\sin(2k\pi x), & \quad x\in (-1,1),\\
     0, &\quad x\in(1,+\infty)\cup (-\infty,-1),
     \end{split}\right.
     \end{equation*}
for $c\in\mathbb C$. Otherwise, for $\lambda=-b^2$, $b>0$,
\begin{equation*}
     g_2(x)=c
     \left\{\begin{split}
      \sin(bx), &\quad x\in (-1,1),\\
     \sin(bx)-\sin(b(x-1)), &\quad x\in(1,+\infty),\\
     \sin(bx)-\sin(b(x+1)), &\quad x\in(-\infty,-1),
     \end{split}\right.
     \end{equation*}
for $d\in\mathbb C$.
This implies that the set of eigenvalues corresponding to odd eigenfunctions is exactly the negative half part of the real line $(-\infty,0]$. All of them are of multiplicity one (restricting to odd eigenfunctions). All these eigenfunctions are merely bounded as $|x|\to\infty$, except when $\lambda=-4k^2\pi^2$, $k=1,2\dots$, that in these cases the eigenfunctions have compact support.

Now we seek  $\lambda$ and $g_1$ even function, solutions of the following eigenvalue problem
\begin{align}\label{eq:even}
\lambda{g_1}= {g_1}_{xx}  -{g_1}_{xx}(0)\mathcal{X}_{(-1,1)}-R{g_1}_{xx}(0)\left[\delta_{-1}+\delta_{1}\right],\quad \lambda = \alpha^2, \quad \alpha\in \C.
\end{align}
As before, we only solve the problem in the half line, with the lateral conditions that $g_1(0)={g_1}_x(0)=0$.

First we study the case $\lambda=0$. One may check that if $R\neq -1$, there are no possible eigenfunctions. However, for $R=-1$ there is a single eigenfunction for the zero eigenvalue, and it is given by:
\begin{equation*}
     g_1(x)=
     \left\{\begin{split}
      c x^2,&\quad x\in (-1,1),\\
      c, &\quad x\in(1,+\infty),\\
      c, &\quad x\in(-\infty,-1).
     \end{split}\right.
     \end{equation*}
Now we assume that $\lambda=\alpha^2$, $\alpha\in\mathbb C\backslash\{0\}$, $\alpha = a + ib$.
Let
\bee
g_1(x)=
\left\{
\begin{split}
&h_1(x)=c_1 e^{\alpha x}+c_2 e^{-\alpha x} + D&\quad\mbox{ in }(0,1),\\
&h_2(x)=d_1 e^{\alpha x}+d_2 e^{-\alpha x} &\quad\mbox{ in }(1,\infty).
\end{split}\right.\eee
The constant $D$ takes into account the term ${g_1}_{xx}(0)\mathcal X_{(-1,1)}$ in the interval $x\in(0,1)$. Substituting $h_1$ into \eqref{eq:even} we get $D = -(c_1+c_2)$. The matching and jump conditions are the following:
\bee
\left\{
\begin{split}
&{h_1}_x(0)=0,\\
&{h_2}_x(1^+) = {h_1}_x(1^-) +R {h_1}_{xx}(0), \\
&h_1(1) = h_2(1).\\
\end{split}\right.\eee
Hence the constants $d_1$, $d_2$, $c_1$, $c_2$ and $\alpha$ satisfy the following constraints:
\bee
\left\{
\begin{split}
&c_1=c_2, \\
&d_1e^\alpha+d_2e^{-\alpha} = c_1e^{\alpha}+ c_1e^{-\alpha} -2c_1, \\
&d_1e^{\alpha} - d_2 e^{-\alpha}  = c_1 e^{\alpha} - c_2 e^{-\alpha} +2 R c_1\alpha.
\end{split}\right.\eee
The eigenfunctions read as
\bee
g_1(x)=
\left\{
\begin{split}
& e^{\alpha x}+ e^{-\alpha x}  -2 &\quad\mbox{ in }(0,1),\\
& e^{\alpha x}( 1-e^{-\alpha} +R\alpha e^{-\alpha}) +e^{-\alpha x}( 1-e^{\alpha} -R\alpha e^{\alpha}) &\quad\mbox{ in }(1,\infty).
\end{split}\right.\eee
Since we are looking for bounded eigenfunctions, we distinguish two cases: $\Re(\alpha)=0$, that is $\lambda =-b^2$ with $b\in(0,\infty)$ and $\Re(\alpha)>0$. As before, the other cases ($b<0$ or $\Re(\alpha)<0$) repeat the same solutions.

If $\Re(\alpha)=0$, $\lambda =-b^2$, $b\in(0,\infty)$, then for any $R\in \R$ the eigenfunction reads as
\bee
g_1(x)=
\left\{
\begin{split}
& \cos(bx)-1 &\quad\mbox{ in }(0,1),\\
& \cos(bx)[1 - \cos(b)+Rb\sin(b)] - \sin(bx)[ \sin(b) + Rb\cos(b)] &\quad\mbox{ in }(1,\infty).
\end{split}\right.\eee
This means that $(-\infty,0)$ is filled again with eigenvalues with even eigenfunctions.

But there are some more eigenvalues with even eigenfunction. If $\Re(\alpha)>0$, then we must impose $1-e^{-\alpha} +R\alpha e^{-\alpha}=0$.
Write $\alpha = a + ib$ with $a>0,b\in \R$, equation $1-e^{-\alpha} +R\alpha e^{-\alpha}=0$ becomes
\be \label{system-1}
\left\{
\begin{split}
& e^{a} \sin(b) +R b=0,  \\
&  e^{a}\cos(b) -1+R a=0.
\end{split}\right.\ee

We study now the two places where instability starts: $R=-1$ and $R=R_0\simeq 9.36$ (it's precise value will be indicated below).

We consider first the case $R=-1$.
Note that for each $R\in (-\infty,-1]$ there exists a unique real positive solution $a$ of the equation $e^a-1+Ra=0$. Then $\lambda=\lambda(R)=a^2\in[0,\infty)$ is an eigenvalue with an even eigenfunction that tends exponentially to zero as $|x|\to\infty$. There are no more real positive eigenvalues.

Now we look at $R=R_0=9.36...$.
We look now at those nonzero $\alpha$ with $a = \pm b$. Such $\alpha$ will correspond to the pure imaginary eigenvalues $\lambda$, with
\bee
\lambda = \left\{
\begin{split}
&  2ia^2  &\quad\mbox{ if } \; a=b,\\
- &   2ia^2   &\quad\mbox{ if } \; a=-b.
\end{split}\right.\eee
Admissible values for $a$ and $R$ are the ones that satisfy the equations
\begin{align}\label{eq:a and M}
\cos(a)-\sin(a) = e^{-a},\quad R = \frac{1-e^a\cos(a)}{a}.
\end{align}
Let us compute now the change of the eigenvalues $\lambda$ across the line $Re(\lambda)=0$ with respect to $R$ (we are still considering the case $Re(\alpha)>0$): via implicit differentiation of equation $e^{\alpha} - 1  +R\alpha=0$ we get
\begin{align*}
\frac{d \lambda}{dR} &=\frac{d \alpha^2}{dR} = -\frac{ 2 \alpha^2}{e^\alpha + R}.
\end{align*}

The rate of change of the real part reads as
\begin{align}\label{derivative}
\Re \left(\frac{d \lambda}{dR} \right)_{|_{a=b}}&= -\frac{ 4 a^2 e^a \sin(a) }{(e^\alpha + R)\overline{(e^\alpha + R)}}=\frac{4a^3R}{|e^\alpha+R|^2},
\end{align}
from the equation $e^{a} \sin(a) +R a=0$. So, we see that if $\lambda$ is purely imaginary, then $\Re \left(\frac{d \lambda}{dM} \right)$ has the same sign as $R$. So, when $R>0$ the eigenvalues can go from the stable to the unstable side, but not conversely, as $R$ increases. And, if $R$ is negative, the other way around. This proves that if the first crossing for $R>0$ happens at $R=R_0$ then there is linear instability for all $R\in(R_0,\infty)$. Also, we will see that the first nonreal crossing for $R<0$ happens at $R_1$ and $R_1<-1$, so we will have linear stability for $-1\le R\le R_0$.

An easy way to analyze the equations (\ref{eq:a and M}) is to look at the graphs of the two functions $R=-e^a\sin(a)/a$ and $R=(1-e^a\cos(a))/a$, and each intersection (except $a=0$) will be a solution. These graphs are represented in the next figure, where the horizontal variable is $a$ and the vertical variable has been taken $\hbox{Argsh}(R)$ instead of $R$, to keep the figure in a reasonable bound.

\centerline{\hbox{\includegraphics[width=10cm,height=6cm]
{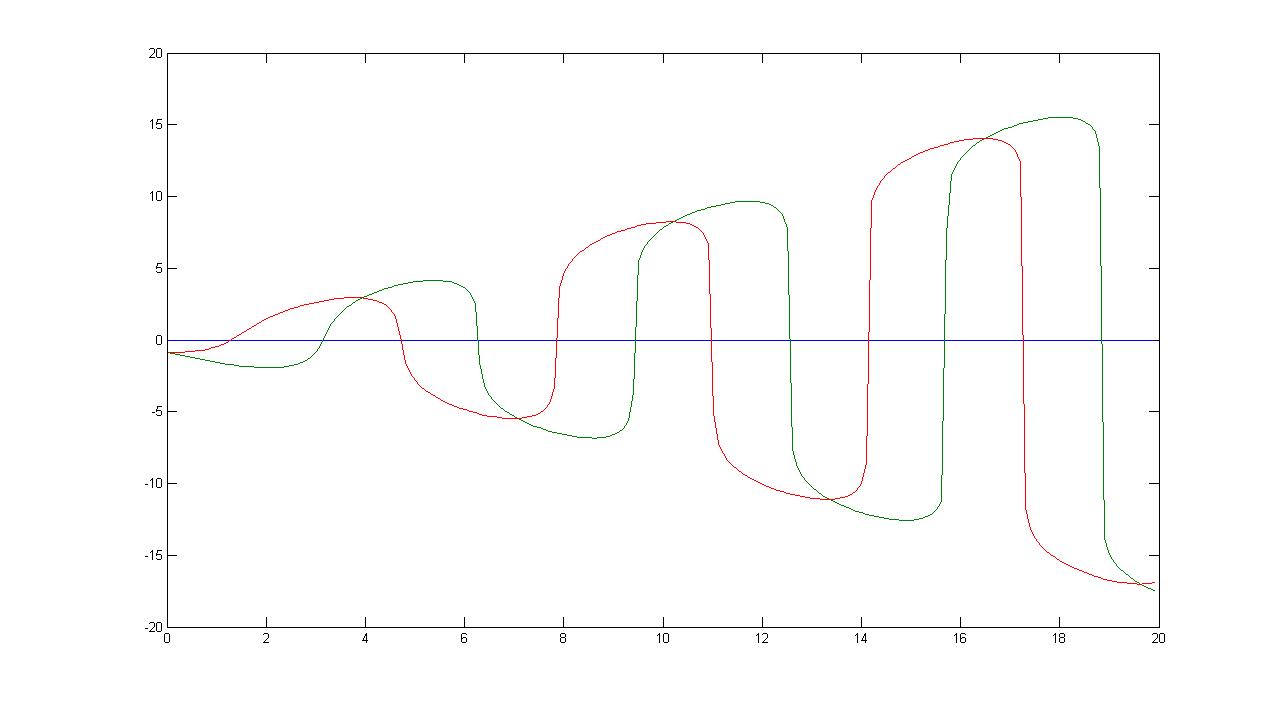}}}

In the figure we see that the first intersections for $R>0$ happens at $a=a_0\simeq 3.94$ and $R=R_0\simeq 9.36$. More exact values are $a_0=3.940733135692915...$ and $R_0= 9.359088829373068...$. For these numbers we obtain $\lambda=\lambda_0\simeq 31.1i$ for the eigenvalue. We also see that there are more (infinitely many more) crossings for larger values of $R$, that would be less significant for the dynamics, but would as well give rise to (unstable) bifurcations of periodic solutions. We also see that he same is true for $R<0$ but the first instabilization happens at a value $R_1\simeq -116$ which is $R_1<-1$.

With this we finish the proof of Theorem \ref{thm1}.

For $R=R_0$ one can calculate the eigenfunction $q_0(s)$ corresponding to $\lambda_0$, and the graph of its real and imaginary parts are represented in the next plot.

\centerline{\hbox{\includegraphics[width=10cm,height=6cm]
{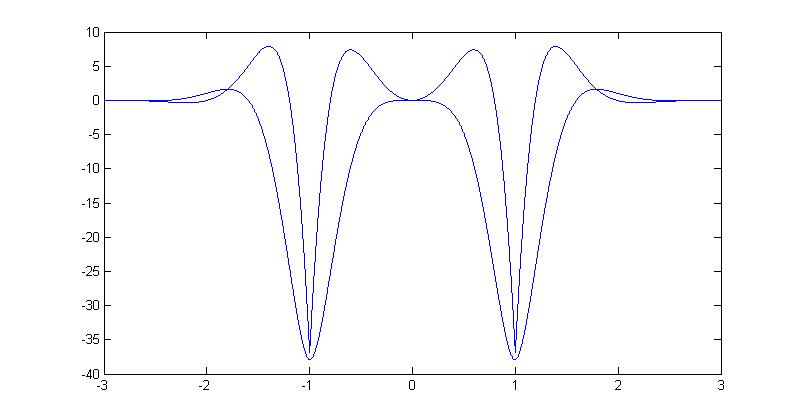}}}

\section{The bifurcating solutions}

This section is devoted to the proof of Theorem \ref{thm2}, and the presentation of the numerical results of Claim \ref{claim}.

Consider again the nonlinear problem (\ref{eqnonlin1}) for the function $w(x,t) = f(x+p(t),t)$.
Let us look at its time-independent solutions for $x\in \R$, that correspond to traveling wave solutions with wave speed $p'(t)=c$ for the original problem. First note that for $p' = 0$, the only equilibrium solutions to \eqref{eqnonlin1} for any constant $R>0$ are the classical trapezoidal shapes (\ref{equilibrium}).

Now we set
$$p' = -\frac{w_{xx}(0)}{w_x(0)}=c,$$
for some constant $c\neq 0$. We seek for continuous functions of the form
\bee
w(x)=
\left\{
\begin{split}
&w_1(x)=a_1+a_2 e^{-c x}&\quad\mbox{ in }(-\infty,-1),\\
&w_2(x)=b_1+b_2 e^{-c x} &\quad\mbox{ in }(-1,1),\\
&w_3(x)=d_1+d_2 e^{-c x}&\quad\mbox{ in }(1,+\infty), \\
\end{split}\right.\eee
that solve the equation
$$
w_{xx} +cw_x(x,t) = w_x(0,t)\left[\delta_{-1}-\delta_{1}\right]+ Rc \left[\delta_{-1}\phi(w(-1))-\delta_{1}\phi(w(1))\right].
$$
Hence $w_1$, $w_2$ and $w_3$ have to satisfy the following matching and jump conditions:
\begin{equation}\label{matching1}
\left\{
\begin{split}
&w_2(0)=0,\\
& w_1(-1^-) = w_2(-1^+),\\
& w_2(1^-) = w_3(1^+),\\
&{w_2}_x(-1^+) - {w_1}_x(-1^-) = {w_2}_x(0) + R c\phi(w(-1)), \\
&{w_3}_x(1^+) - {w_2}_x(1^-) = -{w_2}_x(0) - Rc\phi (w(1)).
\end{split}\right.\end{equation}
First assume that $c>0$. We impose the conditions
\begin{equation}\label{limits1}
w(-\infty)=\rho>0,\quad w(x) \text{ bounded}.
\end{equation}
Then all these these conditions  give that
$$a_1=\rho, \ \ a_2=0, \ \ b_1=\dfrac{\rho}{1-e^c},\ b_2=\dfrac{-\rho}{1-e^c}, \ \ d_1=-\rho\dfrac{\phi(\rho e^{-c})}{\phi(\rho)},\ \ d_2=-\rho+e^c\rho\dfrac{\phi(\rho e^{-c})}{\phi(\rho)}$$
and the relation
\begin{equation}\label{R}
\phi(\rho)=-\rho/R.
\end{equation}

So, the conclusion is that for each value of $c>0$ and each value of $\rho>0$ there exists a unique traveling wave with profile $w^{c,\rho}$ that moves with speed $c$,  $\lim_{x\to-\infty}w^{c,\rho}(x)=\rho$, $\lim_{x\to\infty}w^{c,\rho}(x)= -\rho{\phi(\rho e^{-c})}/{\phi(\rho)}$ and $R$ is given by $R=-\rho/\phi(\rho)$. Note that the profile is constant for $-\infty<x\le -1$ but not for $1\le x<\infty$, and that $\lim_{x\to\infty}w^{c,\rho}(x)$ is not $-\rho$ for most nonlinearities, with the important exception of $\phi=\phi_1$. Note also that as $c\to 0$ with fixed $\rho$ the profile $w^{c,\rho}(x)$ approaches the equilibrium profile $w^\rho$ of (\ref{equilibrium}).

For $c<0$, instead of using $\lim_{x\to-\infty}w(x)=\rho>0$ as a  parameter it is more convenient to parametrize by the limit at $+\infty$, whose value must now be negative, and will be called $-\rho$, where $\rho$ will again be positive. One gets a profile $w^{c,\rho}(x)$ with the property that $w^{c,\rho}(x)=w^{-c,\rho}(-x)$, and we reduce to the previous case.

What we got is that for each pair $(c,\rho)$ with $\rho>0$ there exists a unique traveling wave profile $w^{c,\rho}$ that moves with velocity $c$ and with the following behavior at infinity: if $c>0$ then $\lim_{x\to-\infty} w^{c,\rho}(x)=\rho$, if $c<0$ then $\lim_{x\to\infty } w^{c,\rho}(x)=-\rho$, and if $c=0$ then $\lim_{x\to-\infty }w^{0,\rho}(x)=-\lim_{x\to \infty} w^{0,\rho}(x)=\rho$. This is a two-parameter family solutions that one can say that bifurcates from the equilibria $w^{0,\rho}$, previously called simply $w^\rho$.

For the sake of definiteness, let us describe with detail the profiles one obtains with the first of the nonlinearities defined in (\ref{nonlinearities}), namely $\phi(r)=\phi_1(r)=\hbox{ sign}(r)$. In this case, the traveling wave profiles with $\rho=1$ are given by
\bee
w^{c,1}(x)=
\left\{
\begin{split}
&1 &\quad\mbox{ in }(-\infty,-1),\\
&-\frac{ 1 }{1-e^{c}} (e^{-c x}-1) &\quad\mbox{ in }(-1,1),\\
&-1 - (1-e^c) e^{-c x}&\quad\mbox{ in }(1,+\infty).\\
\end{split}\right.\eee
for $c>0$, and
\bee
w^{c,1}=
\left\{
\begin{split}
&1-(e^{-c}-1)e^{-cx} &\quad\mbox{ in }(-\infty,-1),\\
&-\frac{ 1 }{e^{-c}-1} (e^{-c x}-1) &\quad\mbox{ in }(-1,1),\\
&-1 &\quad\mbox{ in }(1,+\infty).\\
\end{split}\right.\eee
when $c<0$. The graph of the profile $w^{2,1}(x)$ given above is represented in the next figure:

\centerline{\hbox{\includegraphics[width=10cm,height=6cm]
{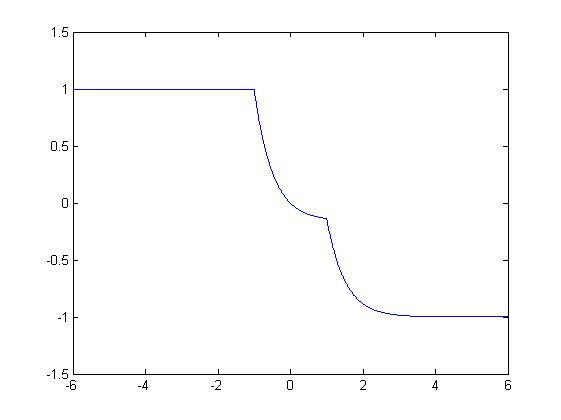}}}

For the general nonlinearity, observe that the formula $R=-\rho/\phi(\rho)$ we have obtained, does not hold for $c=0$. For $c=0$, every solution is an equilibrium, and exists for all values of $R$. But we can take the limit of $-\rho/\phi(\rho)$ when we approach $w^{0,\rho_0}$ from $w^{c,\rho}$. If $\rho_0=1$ and the solution is $w^{0,1}$ (previously called simply $w^1$) we obtain obviously $R=-1/\phi(1)=-1$, in accordance with what was predicted with the linear theory in Theorem \ref{thm1}.

Let us describe with more detail the distribution of the possible values of $R$ for which these solutions exist. For each $R$ these are given by the solutions $\rho>0$ of the equation (\ref{R}). More precisely, for each solution $\rho_R$ of this equation there is a whole family $w^{c,\rho_R}$ of traveling wave profiles, with $c$ moving freely in $c>0$. For this equation to have a solution one needs, first of all, $R<0$, as expected. In the case of the nonlinearity $\phi_1$ for each value of $-\infty<R<0$ there exists a solution of (\ref{R}), namely $\rho_R=-R$. The case of the nonlinearity $\phi_2(r)=r$ is very singular, and solutions of (\ref{R}) only exist for $R=-1$. For $\phi_3(r)=\tanh(r)/\tanh(1)$ a single solution $\rho_R$ exists in the range $-\infty<R<-\tanh(1)\simeq -0.761$. This finishes the proof of Theorem \ref{thm2}.

As we said in the statement, Claim \ref{claim} is merely a well supported conjecture. We do not have a rigorous analytical proof. The difficulties in giving a rigorous proof come from the presence of a continuous part of thew spectrum of the linearized operator that occupies the whole real half-line $(-\infty, 0]$. This does not allow the typical center-manifold techniques for the Hopf bifurcation. Our arguments of support are numerical simulations, made for a problem in a bounded domain.

The numerical simulations have been done for the problem (\ref{eqnonlin1}) but only on the bounded domain $-5<x<5$ and with the boundary conditions $w(-5,t)=1$, $w(5,t)=-1$. We have discretised the spatial domain into small sub-intervals and approximated the differential operators by finite differences. The Dirac's delta functions have been approximated by functions that vanish except at one of the nodes. We have used the Crank-nicolson method to integrate the time evolution.

If one solves the initial value problem with the nonlinearity $\phi=\phi_3$ of (\ref{nonlinearities}) and an initial condition not far from the equilibrium one readily approaches a stable periodic solution for $R$ larger than $R_0$. For $R=12$, for example, one gets a periodic solution of period $T\simeq
0.2088$. The following four plots are the profiles of the solution at times $t=0$, $t=0.0261$, $t=0.0522$ and $t=0.0783$, that cover a half-period.

\centerline{\hbox{\includegraphics[width=7cm,height=5cm]
{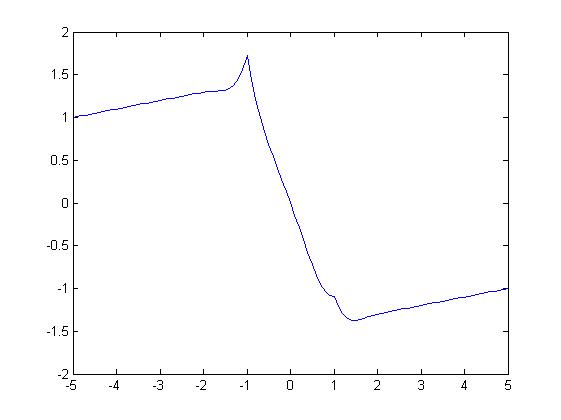}\includegraphics[width=7cm,height=5cm]{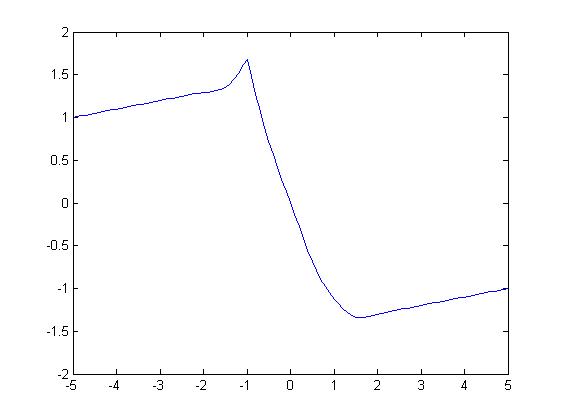}}}
\centerline{\hbox{\includegraphics[width=7cm,height=5cm]
{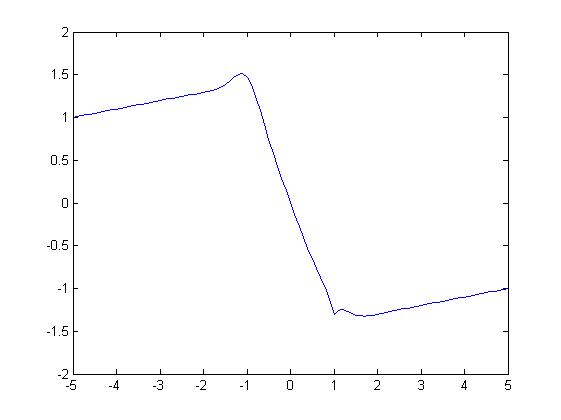}\includegraphics[width=7cm,height=5cm]{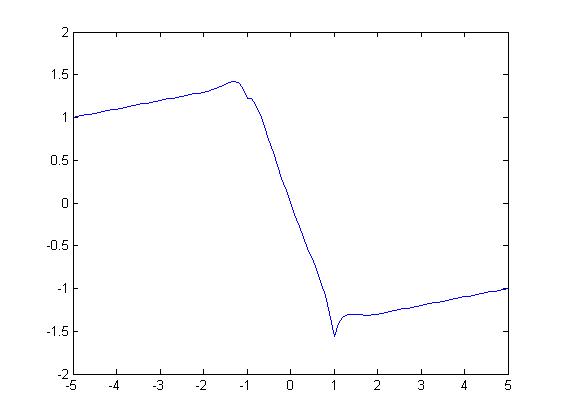}}}

The next half-period, that is $t=0.1044$, $t=0.1305$, $t=0.1566$ and $t=0.1827$ is plotted below, with less resolution.

\centerline{\hbox{
\includegraphics[width=3.5cm,height=2.5cm] {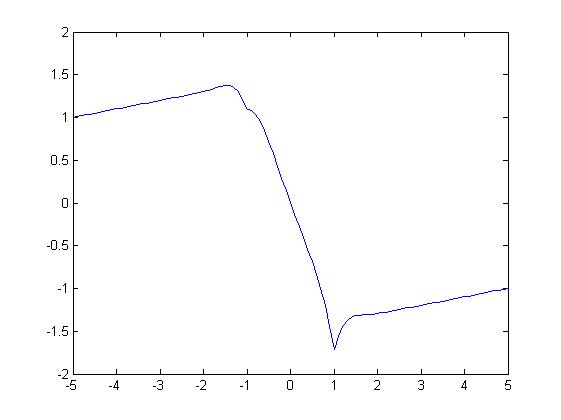}
\includegraphics[width=3.5cm,height=2.5cm] {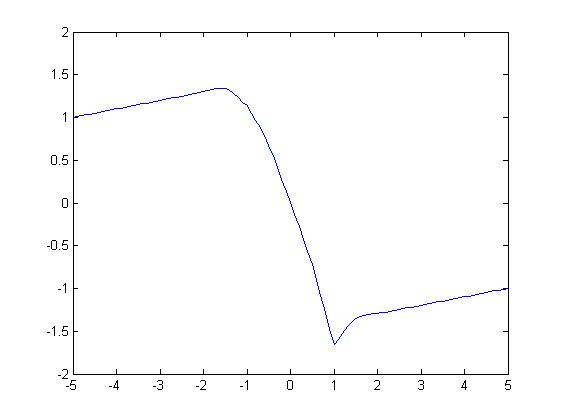}
\includegraphics[width=3.5cm,height=2.5cm] {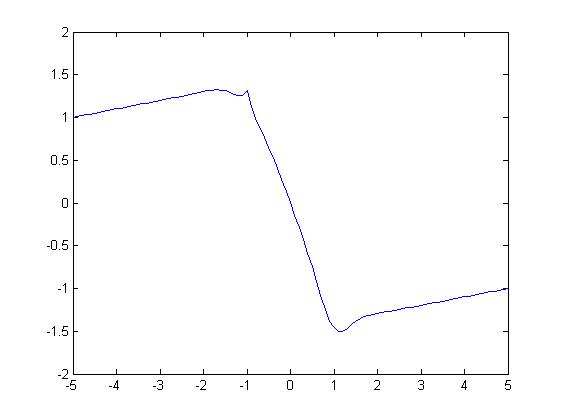}
\includegraphics[width=3.5cm,height=2.5cm] {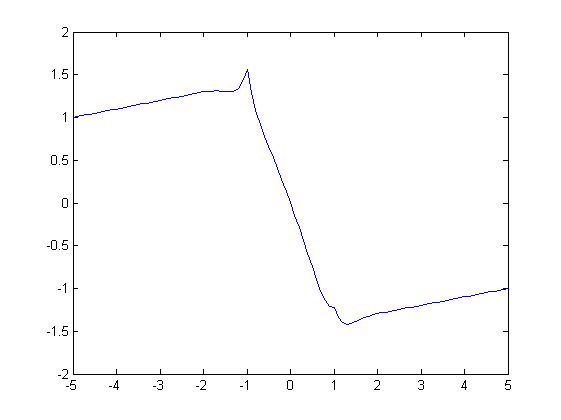}
}}

From the numerical experiments and an analysis of these plots one gets convinced that if $w_p(x,t)$ is this periodic solution then $w_p(x,t+T/2)=-w_p(-x,t)$, a symmetry relation that would deserve an analytical proof. The evolution of the price $p(t)$ along this periodic solution is represented in the following plot:

\centerline{\hbox{\includegraphics[width=7cm,height=5cm]
{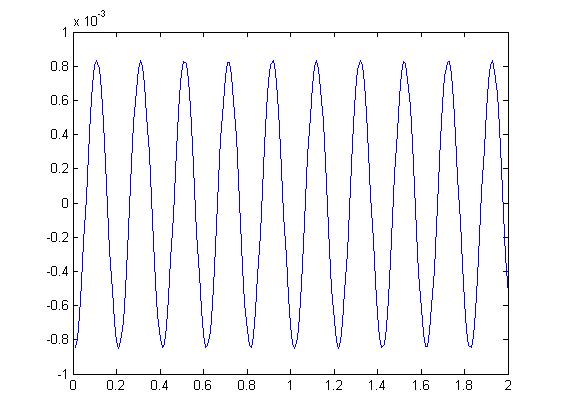}}}

Similar plots are obtained with the same nonlinearity and larger values of $R$.

With the nonlinearity $\phi(r)=\phi_2(r)=r$ the situation is quite different: when the equilibrium becomes unstable, the perturbations seem to grow with time without bound. This is a behavior that would be compatible with a sub-critical Hopf bifurcation.

Finally, with $\phi(r)=\phi_1(r)=\hbox{sign}(r)$ stable oscillations seem to exist but only for for a certain range of values of $R$ after instability. Also, along this range some physically nonsense behaviors do appear, namely the fact that that the periodic profile becomes negative for some negative values of $x$ and positive for some positive values of $x$.

\end{document}